Subject: The Multiplicities of a Dual-thin $Q$-polynomial Association Scheme

Abstract:

Let $Y=(X, \{ R_i \}_{\interval0iD})$ denote a symmetric association scheme,
and assume that $Y$ is $Q$-polynomial with respect to an ordering
$E_0,...,E_D$ of the primitive idempotents.  
In \cite[p.205]{bannai}, Bannai and Ito conjectured that the associated
sequence of multiplicities $m_i$ $(0 \leq i \leq D)$ of $Y$ is unimodal.
We prove that if $Y$ is dual-thin in the sense of Terwilliger,
the sequence of multiplicities satisfies
$m_i \leq m_{i+1}$ and $m_i \leq m_{D-i}$ for $i < D/2$.

----------------------------------------------------------------------

\documentstyle[amssymb,12pt]{article}
\newcommand\newsection[1]{\bigskip\bigskip\refstepcounter{section}
\noindent{\large\bf\thesection\ #1}}
\newcommand\newsubsection[1]{\medskip\refstepcounter{subsection}
{\noindent\bf\thesubsection\ #1.\ }}
\newcommand\newfigur[2]{\medskip\center\refstepcounter{subsection}
{\bf #1 \thesubsection.\ } #2.\ }
\newcommand{\qed}{\rule{1ex}{1ex}}

\newenvironment{theorem}{\newsubsection{Theorem}\sl}{}
\newenvironment{corollary}{\newsubsection{Corollary}\sl}{}
\newenvironment{lemma}{\newsubsection{Lemma}\sl}{}

\newenvironment{conjecture}{\newsubsection{Conjecture}\sl}{}

\newenvironment{proofpart}[1]{\medskip\noindent{\it Proof of {#1}.\ }}
{\mbox{$\qquad\qed$}}

\newcommand\intv[3]{{#1}\leq{#2}\leq{#3}}
\newcommand\interval[3]{{#1}\leq{#2}\leq{#3}}

\newcommand\be{\begin{equation}}
\newcommand\ee{\end{equation}}

\newcommand\bea{\begin{eqnarray}}
\newcommand\eea{\end{eqnarray}}

\newcommand\complexes{\Bbb C}

\begin{document}

\centerline{\bf The Multiplicities of a Dual-thin}
\centerline{\bf $Q$-polynomial Association Scheme}

\bigskip

\centerline{ Bruce E. Sagan }

\centerline{ John S. Caughman, IV }

\bigskip

\centerline{Dept. of Mathematics}

\centerline{Michigan State University }

\centerline{E. Lansing, MI 44824-1027}


\bigskip


\noindent{\bf Abstract.\ }  
Let $Y=(X, \{ R_i \}_{\interval0iD})$ denote a symmetric association scheme,
and assume that $Y$ is $Q$-polynomial with respect to an ordering
$E_0,...,E_D$ of the primitive idempotents.  
In \cite[p.205]{bannai}, Bannai and Ito conjectured that the associated
sequence of multiplicities $m_i$ $(0 \leq i \leq D)$ of $Y$ is unimodal.
We prove that if $Y$ is dual-thin in the sense of Terwilliger,
the sequence of multiplicities satisfies
$m_i \leq m_{i+1}$ and $m_i \leq m_{D-i}$ for $i < D/2$.


\newsection{Introduction}

	For a general introduction to association schemes, we refer
	to \cite{bannai}, \cite{bcn}, \cite{godsil}, or
	\cite{terwSubconstituentI}.  Our notation follows that found in
	\cite{caughmanModules}.

	Throughout this article, $Y=(X, \{ R_i \}_{\interval0iD})$ 
	will denote a symmetric, $D$-class association scheme.  
	Our point of departure is the following well-known result 
	of Taylor and Levingston.

	\begin{theorem} \cite{taylorleving}
	If $Y$ is $P$-polynomial with respect to an ordering 
	$R_0,...,R_D$ of the associate classes,
	then the corresponding sequence of valencies 
		$$k_0, k_1,\ldots, k_D$$
	is unimodal.  Furthermore,
		$$k_i \leq k_{i+1}\quad\mbox{and}\quad k_i \leq  k_{D-i}
		 \quad\mbox{for } i < D/2.\qquad\qed$$
	\end{theorem}

	\noindent Indeed, the sequence is log-concave, as
	is easily derived from the 
	inequalities $b_{i-1} \geq b_i$ and $c_i \leq c_{i+1}$
	$(0 < i < D)$, which are satisfied by the intersection numbers 
	of any $P$-polynomial scheme (cf. \cite[p.\ 199]{godsil}).

	In their book on association schemes, Bannai and Ito made the dual
	conjecture.

	\begin{conjecture} \cite[p.\ 205]{bannai}
	If $Y$ is $Q$-polynomial with respect to an ordering 
	$E_0,...,E_D$ of the primitive idempotents,
	then the corresponding sequence of multiplicites 
		$$m_0,m_1,\ldots,m_D$$
	is unimodal. 
	\end{conjecture}

\smallskip

	\noindent Bannai and Ito further remark that although unimodality of the
	multiplicities follows easily
	whenever the dual intersection numbers satisfy the inequalities 
	$b^*_{i-1} \geq b^*_i$ and $c^*_i \leq c^*_{i+1}$ $(0 < i < D)$,
	unfortunately these inequalities do not always hold. 
	For example, in the Johnson scheme $J(k^2,k)$ we find 
	that $c^*_{k-1} > c^*_k$ whenever $k > 3$.
	However, our main result shows that under a suitable
	restriction on $Y$, the multiplicities satisfy 
	the inequalities
		$$m_i \leq m_{i+1}\quad\mbox{and}\quad m_i \leq  m_{D-i}
		 \quad\mbox{for } i < D/2.$$

	To state our result more precisely, we first review a few definitions.
	Let Mat$_X(\complexes$) denote the $\complexes$-algebra of matrices
	with entries in $\complexes$, where the rows and columns are indexed
	by $X$, and let $A_0$,...,$A_D$ denote the associate matrices 
	for $Y$. Now fix any $x \in X$, and for each integer $i$ 
	$(0 \leq i \leq D)$, let
	$E^*_i = E^*_i(x)$ denote the diagonal matrix in Mat$_X(\complexes$)
	with $yy$ entry
	\begin{equation}\label{e*def}
	(E^*_i)_{yy} = \left\{
	\begin{array}{lr} 1 & \mbox{if } xy \in R_i, \\
			  0 & \mbox{if } xy \not\in R_i.
	\end{array} \right. \; \; \; \; \; \; \; (y \in X).
	\end{equation}
	The {\em Terwilliger algebra} for $Y$ with respect to $x$ is
	the subalgebra $T=T(x)$ of Mat$_X(\complexes$) generated by
	$A_0$,...,$A_D$ and $E^*_0$,...,$E^*_D$.  The Terwilliger 
	algebra was first introduced in \cite{terwSubconstituentI} as 
	an aid to the study of association schemes.
	For any $x \in X$, $T=T(x)$ is a finite dimensional, semisimple 
	$\complexes$-algebra, and is noncommutative in general.
	We refer to \cite{caughmanModules} or \cite{terwSubconstituentI} 
	for more details.  $T$ acts faithfully on the vector space 
	$V:= {\complexes}^X$ by matrix multiplication. 
	$V$ is endowed with the inner product 
	$\langle \; , \; \rangle$ defined by
	$\langle u , v \rangle := u^t \overline{v}$ for all 
	$u, v \in V$.  Since $T$ is semisimple,
	$V$ decomposes into a direct sum of irreducible 
	$T$-modules.

	Let $W$ denote an irreducible $T$-module. 
	Observe that $W = \sum E^*_i W$ (orthogonal direct sum),
	where the sum is taken over all the indices $i$ $(\intv0iD)$ such
	that $E^*_iW \not=0$.
	We set 
		$$d := | \{ i \; : \; E^*_i W \not= 0 \} | -1,$$
	and note that the dimension of $W$ is at least $d+1$. 
	We refer to $d$ as the {\em diameter} of $W$.
	The module $W$ is said to be {\em thin} whenever
	$\dim (E^*_i W) \leq 1$ $(\interval0iD)$.
	Note that $W$ is thin if and only if the diameter of $W$ 
	equals $\dim(W) -1$.
	We say $Y$ is {\em thin} if every irreducible $T(x)$-module is 
	thin for every $x \in X$.

	Similarly, note that $W = \sum E_i W$ (orthogonal direct sum),
	where the sum is  over all $i$ $(\intv0iD)$ such
	that $E_iW \not=0$.
	We define the {\em dual diameter} of $W$ to be
		$$d^* := | \{ i \; : \; E_i W \not= 0 \} | -1,$$
	and note that $\dim W\geq d^*+1$. 
	A {\em dual thin} module $W$ satisfies
	$\dim (E_i W) \leq 1$ $(\interval0iD)$.
	So $W$ is dual thin if and only if  $\dim(W)=d^*+1$.
	Finally, $Y$ is {\em dual thin} if every irreducible $T(x)$-module is 
	dual thin for every vertex $x \in X$.

	Many of the known examples of $Q$-polynomial
	schemes are dual thin.  (See \cite{terwSubconstituentIII} for
	a list.)  Our main theorem is as follows.

\begin{theorem}\label{unim}
	Let $Y$ denote a symmetric
	association scheme which is $Q$-polynomial with respect to 
	an ordering $E_0,...,E_D$ of the primitive idempotents.  If $Y$ is 
	dual-thin, then the multiplicities satisfy 
		$$m_i \leq m_{i+1}\quad\mbox{and}\quad m_i \leq  m_{D-i}
		 \quad\mbox{for } i < D/2.$$
\end{theorem}

	\noindent The proof of Theorem \ref{unim} is contained in the
	next section. 

	We remark that if $Y$ is bipartite $P$- and $Q$-polynomial, then
	it must be dual-thin and
	$m_i = m_{D-i}$ for $i< D/2$.  So Theorem \ref{unim} implies
	the following corollary. 
	(cf. \cite[Theorem 9.6]{caughmanSpectra}).
\begin{corollary}
	Let $Y$ denote a symmetric
	association scheme which is bipartite $P$- and $Q$-polynomial
	with respect to an ordering $E_0,...,E_D$ of the
	primitive idempotents. 	Then the corresponding sequence of
	multiplicites  
		$$m_0,m_1,\ldots,m_D$$
	is unimodal.\qquad\qed
\end{corollary}
 

\newsection{Proof of the Theorem}

	Let $Y=(X, \{ R_i \}_{\interval0iD})$ denote a symmetric 
	association scheme which
	is $Q$-polynomial with respect to the ordering $E_0,...,E_D$
	of the primitive idempotents.  Fix any $x \in X$ and 
	let $T=T(x)$ denote the Terwilliger algebra for $Y$ with
	respect to $x$.  Let $W$ denote any
	irreducible $T$-module.  We define the {\em dual endpoint } of 
	$W$ to be the integer $t$ given by
\begin{eqnarray}
	t &:=& \min \{ i \; : \; 0 \leq i \leq D, \; E_i W \not= 0 \}. 
		\label{dep} 
\end{eqnarray}
	We observe that $0 \leq t \leq D-d^*, \;$ where $d^*$ denotes 
	the dual diameter of $W$.

\begin{lemma}\label{qpolythin}\cite[p.385]{terwSubconstituentI}
	Let $Y$ be a symmetric association scheme which
	is $Q$-polynomial with respect to the ordering $E_0,...,E_D$
	of the primitive idempotents. Fix any $x \in X$, and write 
	$E^*_i = E^*_i(x)$ $(\interval0iD)$, $T=T(x)$.
	Let $W$ denote an irreducible 
	$T$-module with dual endpoint $t$.  Then
\begin{enumerate}
\begin{item} $E_i W \not= 0 \hspace{2em} \mbox{ iff } \hspace{2em}
	\interval{t}i{t+d^*} 
	\hspace{3em} (\interval0iD)$.
\end{item}
\begin{item} Suppose $W$ is dual-thin.  Then $W$ is thin, and $d = d^*$.\qquad\qed
\end{item}
\end{enumerate}
\end{lemma}

\begin{lemma}\cite[Lemma 4.1]{caughmanModules}\label{2t+dlem1}
	Under the assumptions of the previous lemma, the dual
	endpoint $t$ and diameter $d$ of any irreducible $T$-module satisfy
$$
	2t+d \geq D.\qquad\qed
$$
\end{lemma}

\begin{proofpart}{Theorem \ref{unim}}
Fix any $x \in X$, and let $T=T(x)$ denote the Terwilliger algebra for $Y$
with respect to $x$.  Since $T$ is semisimple, there exists a positive
integer $s$ and irreducible $T$-modules $W_1$, $W_2$,...,$W_s$ such that
\begin{equation}
	V = W_1 + W_2 + \cdots + W_s \hspace{2.5em} 
	\mbox{(orthogonal direct sum).}
\end{equation}
For each integer $j$, $1 \leq j \leq s$, let $t_j$ (respectively,
$d_j^*$) denote the dual endpoint (respectively, dual diameter) of $W_j$.
Now fix any nonnegative integer $i< D/2$.  Then for any $j$,
$1 \leq j \leq s$, 
$$
\begin{array}{rcll}
E_i W_j \not= 0 &\Rightarrow& t_j \leq i 
			&\mbox{(by Lemma~\ref{qpolythin}(i))} \\
		&\Rightarrow& t_j<i+1\leq D-i \leq D- t_j 
			&\mbox{(since $i<D/2$)} \\
		&\Rightarrow& t_j<i+1\leq D-i\leq t_j + d^*_j
			&\mbox{(by Lemmas~\ref{qpolythin}(ii),
			\ref{2t+dlem1}}) \\
		&\Rightarrow& E_{i+1} W_j \not= 0 \mbox{ and }
			E_{D-i} W_j \not= 0 
			&\mbox{(by Lemma~\ref{qpolythin}(i)).}
\end{array}
$$
So we can now argue that, since $Y$ is dual thin,
\begin{eqnarray*}
\dim (E_i V) &=& | \{ j\; :\; 0 \leq j \leq s, \, E_i W_j \not= 0 \} | \\
		&\leq& | \{ j\; :\ ;0 \leq j \leq s, \, E_{i+1} W_j \not= 0 \} | \\
		&=& \dim (E_{i+1} V).
\end{eqnarray*}
In other words, $m_i \leq m_{i+1}$.  Similarly, 
\begin{eqnarray*}
\dim (E_i V) &=& | \{ j\; :\; 0 \leq j \leq s, \, E_i W_j \not= 0 \} | \\
		&\leq& | \{ j\; :\; 0 \leq j \leq s, \, E_{D-i} W_j \not= 0 \} | 	\\
		&=& \dim (E_{D-i} V) 
\end{eqnarray*}
This yields $m_i \leq m_{D-i}$. 
\end{proofpart}

\end{document}